\documentclass[12pt]{article}
\usepackage[T1]{fontenc}
\usepackage[francais]{babel}
\usepackage{graphicx}
\DeclareGraphicsExtensions{.jpg}
\makeatletter
\def\@begintheorem#1#2{\trivlist
      \item[\hskip \labelsep{\bf #1\ #2}]\rm}
\def\@opargbegintheorem#1#2#3{\trivlist
      \item[\hskip \labelsep{\bf #1\ #2\ (#3)}]\rm}
\makeatother         
\title{La m\'ethode d'it\'eration rationnelle de Georges Lema\^itre (document de travail)\\
The rational iteration method by Georges Lema\^itre (working papers)}
\author{Herv\'e  Le Ferrand \footnote{Institut de Math\'ematiques de Bourgogne, Universit\'e de Bourgogne, France, leferran@u-bourgogne.fr}}
\date{\today}

\begin{document}
\maketitle

{\bf Abstract}\\
In 1942, the famous astronomer and physicist Belgian Georges Lema\^itre (1894-1966) published in the Bulletin of the class of Sciences of the Royal Academy of Belgium, an article in which he proposes a method that he calls "rational iteration". He rediscovers both the $\Delta^{2}$ Aitken process and the Steffensen's method.

We will see how the scientific concerns of Georges Lema\^itre in the year 1942 (resolutions of differential equations) have led to this rediscovery.

{\bf Key words}
Lema\^itre, Aitken, Steffensen, $\Delta^{2}$, fixed point, ordinary differential equation, Picard iterative process.

\vspace{1cm}

{\bf R\'esum\'e}\\
En 1942, le c\'el\`ebre astronome et physicien belge Georges Lema\^itre (1894-1966) publie dans le Bulletin de la classe des sciences de l'Acad\'emie Royale de Belgique, un article \cite{Lem1} dans lequel il formule une m\'ethode qu'il nomme {\it it\'eration rationnelle}. Il red\'ecouvre \`a la fois le processus du {\it $\Delta^{2}$ de Aitken} et la {\it m\'ethode de Steffensen}. C'est en 1927 que le math\'ematicien n\'eo-z\'elandais Alexander Craig Aitken (1895-1967) publie son fameux proc\'ed\'e d'acc\'el\'eration de convergence \cite{Aitken1}. Quelques ann\'ees plus tard, en 1933, le math\'ematicien danois Johan Frederik Steffensen (1873-1961) expose sa m\'ethode \cite{Stef1}. La m\'ethode de Steffensen est vue d'ailleurs \`a pr\'esent comme un corollaire du $\Delta^{2}$ de Aitken\footnote{ Steffensen ne mentionne pas Aitken dans son article. Le compte-rendu paru dans Zentralblatt sur l'article de 1933 ne fait pas non plus mention de Aitken.}. 

Si le compte-rendu de l'article de Georges Lema\^itre paru dans le Jahrbuch de 1942 mentionne la m\'ethode de Steffensen\footnote{Le compte-rendu de Zentralblatt est similaire. Dans les deux cas, il n'est pas question du proc\'ed\'e d'Aitken.}, Lema\^itre ne fait pas r\'ef\'erence aux deux math\'ematiciens. Cela peut para\^itre surprenant, mais nous allons voir que les conditions de vie et de travail de Georges Lema\^itre en cette ann\'ee de guerre 1942 ainsi que ses pr\'eoccupations scientifiques expliquent cela.

{\bf Mots cl\'es}
Lema\^itre, Aitken, Steffensen, $\Delta^{2}$, point fixe, \'equation diff\'erentielle ordinaire, m\'ethode it\'erative de Picard.
\newpage

\tableofcontents
\newpage

\section{Le $\Delta^{2}$ de Aitken et la m\'ethode de Steffensen}

\subsection{L'article de Aitken}
Le r\'esum\'e de l'article de Aitken\footnote{Le r\'esum\'e est disponible sur :  http://dx.doi.org/10.1017/S0370164600022070 } est le suivant :

\begin{quotation}
\begin{it}
Proceedings\\
XXV-On Bernoulli's Numerical Solution of Algebraic Equations\\
A. C. Aitken

The aim of the present paper is to extend Daniel Bernoulli's method of approximating to the numerically greatest root of an algebraic equation. On the basis of the extension here given it now becomes possible to make Bernoulli's method a means of evaluating not merely the greatest root, but all the roots of an equation, whether real, complex, or repeated, by an arithmetical process well adapted to mechanical computation, and without any preliminary determination of the nature or position of the roots. In particular, the evaluation of complex roots is extremely simple, whatever the number of pairs of such roots. There is also a way of deriving from a sequence of approximations to a root successive sequences of ever-increasing rapidity of convergence.
\end{it}
\end{quotation}

Ainsi le point de d\'epart du travail de Aitken est le m\'emoire de Daniel Bernoulli (1700-1782) {\it Observationes de seriebus quae formantur ex additione vel subtractione quacunque terminorum se mutuo consequentium}\cite{Bernoulli}. A une \'equation alg\'ebrique, mise sous une forme particuli\`ere, Bernoulli fait correspondre une suite r\'ecurrente lin\'eaire. En utilisant une \'ecriture moderne, si l'\'equation est $1=a_{1}x+\cdots+ a_{p}x^{p}$, la suite associ\'ee est $u_{n+p}=a_{1}u_{n+p-1}+\cdot+a_{p}u_{n}$ ($n\geq 0$)\footnote{C'est d'ailleurs la pr\'esentation choisie par Henrici\cite{Henrici}}. Avec de bonnes hypoth\`eses, le rapport $\displaystyle{\frac{u_{n}}{u_{n+1}}}$ tend vers la solution de plus petit module de l'\'equation. Une analyse de cette m\'ethode est donn\'ee dans \cite{Chabert}, pages 254-259. Les auteurs indiquent que Leonhard Euler dans \cite{Euler} donne une preuve de la m\'ethode de Bernoulli en consid\'erant une fraction rationnelle de d\'enominateur $1-a_{1}x-\cdots -a_{p}x^{p}$ qu'il d\'eveloppe en s\'erie enti\`ere de deux fa\c cons. Tout d'abord, il fait remarquer que les coefficients du d\'eveloppement suivent la relation de r\'ecurrence\footnote{\`a partir d'un certain rang.}. Puis il obtient le d\'eveloppement en d\'ecomposant la fraction rationnelle en \'el\'ements simples en faisant donc  appara\^itre les racines de l'\'equation. Euler d'\'ecrit d'ailleurs le cas d'une \'equation \`a racines simples : l'\'el\'ement simple $\displaystyle{\frac{1}{1-pz}}$ se d\'eveloppe sous la forme$\sum_{k\geq 0}p^{k}z^{k}$, ce qui permet d'avoir la r\`egle de Bernoulli \footnote{Euler indique de plus que pour obtenir la racine de plus grand module de $P(z)$, on applique la m\'ethode au polyn\^ome $z^{p}P(\frac{1}{z})$, $P$ \'etant de degr\'e $p$.}. On retrouve toute cette d\'emarche expos\'ee dans \cite{Henrici}, pages 584-588.

Aitken\footnote{Une analyse de son article est faite dans \cite{Chabert} pages 492-497.} commence par g\'en\'eraliser la m\'ethode pr\'ec\'edente  de Daniel Bernoulli dite {\it m\'ethode des s\'eries r\'ecurrentes}. Partant d'une \'equation $z^{p}+b_{p-1}z^{p-1}+\cdots+b_{1}z+b_{0}=0$ et de la suite r\'ecurrente $s_{n+p}+b_{p-1}s_{n+p-1}+\cdots+b_{1}s_{n+1}+b_{0}s_{n}=0$, il consid\`ere les rapports :
\begin{equation}
\frac{
\left\vert
\begin{array}{lrrl}
s_{n+1}&s_{n+2}&\cdots&s_{n+m}\\
s_{n}&s_{n+1}&\cdots&s_{n+m-1}\\
\cdots&\cdots&\cdots&\cdots \\
s_{n-m+2}&s_{n-m+3}&\cdots&s_{n+1}
\end{array}\right\vert
}
{
\left\vert
\begin{array}{lrrl}
s_{n}&s_{n+1}&\cdots&s_{n+m-1}\\
s_{n-1}&s_{n}&\cdots&s_{n+m-2}\\
\cdots&\cdots&\cdots&\cdots \\
s_{n-m+1}&s_{n-m+2}&\cdots&s_{n}
\end{array}\right\vert
}
\end{equation}
Aitken prouve que la limite de ce rapport quand $n\to +\infty$ est $z_{1}z_{2}\cdots z_{m}$ o\`u $z_{1}, z_{2},\ldots, z_{m}$ sont les $m$ premi\`eres racines de l'\'equation, suppos\'ees telles que $\vert z_{1}\vert>\vert z_{2}\vert>\cdots>\vert z_{m}\vert$.

En permutant les lignes dans chacun des d\'eterminants, on obtient le rapport :
\begin{equation}
\frac{
\left\vert
\begin{array}{lrrl}
s_{n-m+2}&s_{n-m+3}&\cdots&s_{n+1}\\
\cdots&\cdots&\cdots&\cdots \\
s_{n}&s_{n+1}&\cdots&s_{n+m-1}\\
s_{n+1}&s_{n+2}&\cdots&s_{n+m}
\end{array}\right\vert
}
{
\left\vert
\begin{array}{lrrl}
s_{n-m+1}&s_{n-m+2}&\cdots&s_{n}\\
\cdots&\cdots&\cdots&\cdots \\
s_{n-1}&s_{n}&\cdots&s_{n+m-2}\\
s_{n}&s_{n+1}&\cdots&s_{n+m-1}
\end{array}\right\vert
}
\end{equation}
qui est en fait le quotient de deux d\'eterminants de Hankel (1839-1873) associ\'es \`a la s\'erie $\sum_{n\geq 0}s_{n}z^{n}$ : $\displaystyle{\frac{H_{m}^{n-m+2}}{H_{m}^{n-m+1}}}$.
Rappelons que l'expression g\'en\'erale des d\'eterminants de Hankel est :
\begin{equation}
H_{j}^{i}=\left\vert
\begin{array}{lrrl}
s_{i}&s_{i+1}&\cdots&s_{i+j-1}\\
s_{i+1}&s_{i+2}&\cdots&s_{i+j}\\
\cdots&\cdots&\cdots&\cdots \\
s_{i+j-1}&s_{i+j}&\cdots&s_{i+2j-2}
\end{array}\right\vert.
\end{equation}
Peter Henrici donne une d\'emonstration de ce r\'esultat de convergence dans \cite{Henrici}, pages 596-602 \footnote{Il fait d'ailleurs r\'ef\'erence \`a l'article de Aitken.}. Il lie ensuite ce r\'esultat \`a l'algorithme {\it quotient-diff\'erence} de Heinz Rutishauser (1918-1970)\cite{Rutis1}.

\subsection{Acc\'el\'eration de la convergence}

Aitken cherche ensuite \`a am\'eliorer la vitesse de convergence de suites qu'il a pr\'ealablement consid\'er\'ees. Ce qu'il fait en consid\'erant les diff\'erences finies des premier et second ordre de la suite $\displaystyle{\left(t_{n}=\frac{s_{n+1}}{s_{n}}\right)}$. En utilisant son r\'esultat ci-dessus pour $m=2$, il montre que les diff\'erences premi\`eres se comportent comme les termes d'une suite g\'eom\'etrique de raison $\displaystyle{\frac{z_{2}}{z_{1}}}$. Aitken pose alors la relation :
\begin{equation}
\frac{z_{1}-t_{n+2}}{z_{1}-t_{n+1}}=\frac{\Delta t_{n+1}}{\Delta t_{n}}
\end{equation}
En r\'esolvant par rapport \`a $z_{1}$, il vient\footnote{Dans l'article de Aitken, le num\'erateur de la fraction qui suit est \'ecrit sous la forme d'un d\'eterminant $2\times 2$.} :
\begin{equation}
z_{1}=\frac{t_{n}t_{n+2}-t_{n+1}^{2}}{\Delta^{2}t_{n}}.
\end{equation}
Le terme de droite est ce que l'on nomme aujourd'hui le transform\'e par le $\Delta^{2}$ d'Aitken de la suite $(t_{n})$. On esp\`ere que cette nouvelle suite va converger vers $z_{1}$ plus vite que la suite initiale. Aitken propose ensuite de faire subir \`a cette nouvelle suite la transformation du $\Delta^{2}$ pour avoir une convergence plus rapide. Comme le rappelle Claude Brezinski dans \cite{Brez3}, Aitken indique que la transformation du $\Delta^{2}$ a d\'ej\`a \'et\'e consid\'er\'ee par d'autres auteurs mais pas dans un but d'acc\'el\'eration de la convergence. On peut nuancer un peu cela, toujours dans \cite{Brez3} (voir aussi \cite{Brez1} et \cite{Osada}), Claude Brezinski indique que le math\'ematicien japonais Takakazu Seki (?-1708) introduit la transformation du $\Delta^{2}$ vers 1680 dans le calcul du volume d'une sph\`ere. D'apr\`es Osada, la m\'ethode de Seki correspond \`a une acc\'el\'eration de convergence\footnote{Osada \'ecrit : {\it \og From the numerical analysis view point, the formula  is the m-panels trapezoidal rule for the numerical integration $\displaystyle{\int_{0}^{D}4y(D-y)dy=\frac{2}{3}D^{3}}$. Therefore we can say that Seki accelerated the trapezoidal rule by the Aitken
$\Delta^{2}$ process.\fg}}.

\subsection{La m\'ethode de Steffensen}

L'article de Steffensen \cite{Stef1} est publi\'e en 1933 dans les Skandinavisk Aktuarietidskrift, une revue destin\'ee aux actuaires\footnote{Steffensen a d\'ebut\'e sa carri\`ere dans les assurances \`a Copenhague.} dont il a \'et\'e un des \'editeurs pendant trente ans\cite{Ogborn}. Contrairement \`a Aitken et \`a Lema\^itre, le point de d\'epart de Steffensen n'est pas un m\'emoire original d'un math\'ematicien c\'el\`ebre. Steffensen fait r\'ef\'erence \`a deux ouvrages : {\it The calculus of observations} de Whittaker-Robinson\cite{Whit} et {\it Vorlesungen \"uber numerisches Rechnen } de Runge-K\"onig. Steffensen consid\`ere une \'equation scalaire $f(x)=x$ dont il cherche \`a approcher une racine r\'eelle. Il d\'ebute son article par quelques consid\'erations sur la suite des it\'er\'ees $x_{n+1}=f(x_{n})$ selon si $f$ est croissante ou d\'ecroissante (sur un intervalle). Ceci lui permet d'encadrer la racine. Steffensen introduit la notion de {\it suites it\'er\'ees oppos\'ees}. 

Ensuite Steffensen raisonne en termes d'interpolation polyn\^omiale. Il interpole la fonction $f$ aux points (it\'er\'es) $x_{0}$, $x_{1}=f(x_{1})$ et $x_{2}=f(x_{1})$\footnote{Comme l'indique Steffensen, tout choix de trois points cons\'ecutifs convient pour l'\'etude.} en utilisant les diff\'erences divis\'ees\footnote{voir par exemple \cite{Davis} pages 39-40 et 64-65.} :
\begin{equation}
x=f(x)=f(x_{0})+(x-x_{0})\left\lbrack f(x_{0}),f(x_{1})\right\rbrack+(x-x_{0})(x-x_{1})\left\lbrack f(x),f(x_{0}),f(x_{1})\right\rbrack
\end{equation}
ou encore
\begin{equation}
x\approx  x_{1}+(x-x_{0})\frac{x_{2}-x_{1}}{x_{1}-x_{0}}
\end{equation}
Cette derni\`ere expression conduit \`a :
\begin{equation}
x\approx x_{0}-\frac{(\Delta x_{0})^{2}}{\Delta^{2}x_{0}}.
\end{equation}
Steffensen d\'eveloppe alors plusieurs exemples num\'eriques reposant sur l'algorithme dont nous donnons une \'ecriture classique :
\begin{equation}
x_{0},\ \mbox{pour}\ k=0\ \mbox{\`a}\ n\ \mbox{faire}\ y_{0}:=x_{k},\ y_{1}:=f(y_{0}),\ y_{2}:=f(y_{1}),\ x_{k+1}:=y_{0}-\frac{(\Delta y_{0})^{2}}{\Delta^{2}y_{0}}.
\end{equation}
Sur les exemples num\'eriques propos\'es par Steffensen, l'am\'elioration de la vitesse de convergence est bien visible.

Si deux exemples propos\'es par Steffensen portent sur des calculs d'int\'er\^ets, Steffensen devait \^etre aussi motiv\'e par des questions de {\it Dynamique de populations}. En effet, en 1930, Steffensen reprend une question sur le {\it probl\`eme d'extinction} \'etudi\'ee alors par Agner Krarup Erlang (1878-1929) \cite{Baca}. On revient toujours \`a une recherche de point fixe.

\section{L'it\'eration rationnelle de Georges Lema\^itre}
\subsection{L'ann\'ee 1942}
Quel est le contexte historique et quelles sont les pr\'eoccupations scientifiques de Georges Lema\^itre en 1942 ? Nous nous appuyons ici essentiellement sur la biographie de Dominique Lambert 
\cite{Lambert}. En Mai 1940, la Belgique est envahie par l'arm\'ee allemande. La biblioth\`eque de l'Universit\'e Catholique de Louvain est d\'etruite par le feu comme en Ao\^ut 1914\footnote{site de l'Universit\'e Catholique de Louvain}. Apr\`es avoir tent\'e de gagner l'Angleterre avec plusieurs membres de sa famille, Georges Lema\^itre revient en 
Belgique en Juin 1940. Il d\'em\'enage alors une s\'erie de machines \`a calculs m\'ecaniques de type Mercedes\footnote{Lambert op.cit., page 218.}. Georges Lema\^itre utilise ces machines pour ses calculs et veut \'eviter de plus qu'elles ne tombent entre les mains de l'occupant. Georges Lema\^itre est sensibilis\'e au moment de ses \'etudes universitaires \`a cette question du calcul sur machine par ses professeurs (puis coll\`egues), Edouard Goedseels (1857-1928) et surtout Maurice Alliaume (1882-1931)\footnote{Lambert op.cit., pages 42-43.}. A partir des ann\'ees 40, Georges Lema\^itre va consacrer une part importante de son activit\'e scientifique au calcul sur machines. Le point d'orgue est l'introduction -et la programmation- du premier ordinateur \`a 
l'Universit\'e  Catholique de Louvain en 1958\footnote{voir aussi la biographie de Georges Lema\^itre sur le site de l'UCL.}. Dans \cite{Jones}, Pierre-Jacques Courtois dans sa contribution {\it The Belgian Electronic Mathematical Machine (1951-1962) : An Account}, indique que Charles Manneback (1894-1975) a port\'e le projet de construction d'un ordinateur belge au d\'ebut des ann\'ees 50. Or Manneback, professeur de Physique-Math\'ematique et d'Electricit\'e \`a l'\'ecole d'ing\'enieurs de Louvain, est un coll\`egue et ami de Georges Lema\^itre. Dans les ann\'ees cinquante, en Suisse, un projet du m\^eme type voit le jour. Le math\'ematicien Heinz Rutishauser, cit\'e plus haut et qui a g\'en\'eralis\'e -et impl\'ement\'e- la m\'ethode de Aitken, participe \`a la cr\'eation du premier ordinateur suisse, ERMETH \cite{Bauer}. 

\subsection{La m\'ethode}
Georges Lema\^itre propose une {\it modification} de la m\'ethode des approximations successives, $x_{n+1}=f(x_{n})$ pour r\'esoudre l'\'equation $x=f(x)$\footnote{Lema\^itre indique que cette m\'ethode a \'et\'e \'etendue \`a la r\'esolution d'\'equations diff\'erentielles par Emile Picard (1856-1941). C'est justement la motivation de Lema\^itre comme on le voit d'ailleurs \`a la fin de son article. } Notons que Lema\^itre fait r\'ef\'erence comme Steffensen aux deux ouvrages : {\it The calculus of observations} de Whittaker-Robinson\cite{Whit} et {\it Vorlesungen \"uber numerisches Rechnen } de Runge-K\"onig. Deux ouvrages classiques de l'\'epoque, qui devaient se trouver \^etre facilement accessibles. A l'instar de Aitken, Georges Lema\^itre construit son propos \`a partir des travaux d'un math\'ematicien illustre, en l'occurence Carl Friedrich Gauss (1777-1855). Ceci n'est pas surprenant. Dominique Lambert \footnote{Lambert op.cit, page 222.} souligne que la lecture dans les textes originaux (en latin) faisait partie des habitudes de Lema\^itre, habitude h\'erit\'ee de Henri Bosmans (1852-1928)\footnote{Concernant la vie et l'oeuvre du p\`ere j\'esuite Henri Bosmans, on peut se r\'ef\'erer \`a \cite{Hermans}.}. D'ailleurs Lambert pr\'ecise que Georges Lema\^itre {\it priv\'e de ses contacts internationaux} durant la Seconde Guerre mondiale, s'est plong\'e en particulier dans la {\it M\'ecanique Analytique} de Lagrange (1736-1813) et dans {\it Les m\'ethodes nouvelles de la m\'ecanique c\'eleste} de Henri Poincar\'e (1854-1912). 

Georges Lema\^itre adapte \`a l'it\'eration $x=f(x)$ la m\'ethode de Gauss. Gauss s'int\'eresse \`a un syst\`eme 
\begin{equation}
\left\lbrace
\begin{array}{c}
X(x,y)=0\\
Y(x,y)=0
\end{array}\right.
\end{equation}
dont il connait trois valeurs approch\'ees $(x_{n},y_{n})$, $n=1,2,3$. Il s'agit de d\'eterminer \`a partir des ces trois valeurs, une valeur pour $(x,y)$. Gauss lin\'earise le probl\`eme en consid\'erant les d\'eveloppements de Taylor \`a l'ordre $1$, sans reste, au point $(x,y)$ des fonctions $X$ et $Y$. En \'evaluant ces deux d\'eveloppements en $(x_{n},y_{n})$ pour $n=1,2,3$, on obtient trois syst\`emes lin\'eaires d'inconnues respectives $(x-x_{n},y-y_{n})$ ($n=1,2,3$). Gauss donne ensuite l'expression pour un $n$ fix\'e, de $x-x_{n}$ en fonction de $X_{n}=X(x_{n},y_{n})$, de $Y_{n}=Y(x_{n},y_{n})$ et des coefficients des d\'eveloppements de Taylor. Ceux-ci \'etant inconnus, Gauss \'ecrit une {\it condition de compatibilit\'e} entre les trois expressions sous la forme d'un d\'eterminant nul. Georges Lema\^itre reprend cette d\'emarche. Il lin\'earise $f$,
\begin{equation}
\left\lbrace
\begin{array}{l}
x_{2}=f(x_{1})=\alpha +\beta x_{1}\\
x_{3}=f(x_{2})=\alpha +\beta x_{2}\\
x=f(x)=\alpha +\beta x
\end{array}\right.
\end{equation}
puis donne  la condition de compatibilit\'e du syst\`eme :
\begin{equation}
\left\vert
\begin{array}{ccc}
x_{2}&1&x_{1}\\
x_{3}&1&x_{2}\\
x&1&x
\end{array}\right\vert=0.
\end{equation}
En d\'eveloppant le d\'eterminant selon la derni\`ere ligne, on obtient :
\begin{equation}
x=\frac{x_{1}x_{3}-x_{2}^{2}}{x_{1}-2x_{2}+x_{3}}
\end{equation}
c'est \`a dire la transformation du $\Delta^{2}$ d'Aitken !

Si Georges Lema\^itre ne donne pas d'exemple num\'erique, il consid\`ere bien la m\'ethode comme une m\'ethode d'acc\'el\'eration de la convergence. En effet, il conclut :
\begin{quotation}
\begin{it}
On peut donc employer la formule d'it\'eration rationnelle pour acc\'el\'erer la convergence l\`a o\`u elle existe.
\end{it}
\end{quotation}
Cette phrase suit les consid\'erations faites par  Georges Lema\^itre sur la r\'esolution d'une \'equation diff\'erentielle\footnote{probl\`eme avec condition initiale.} $\displaystyle{\frac{dx}{dt}=f(x,t)}$ pour laquelle on conna\^it  une solution approch\'ee $x_{2}(t)$. En suivant le processus it\'eratif de Picard\cite{Picard}, pages 197-200
\begin{equation}
x_{3}=\int f(x_{2},t)dt\ ;\ f(x_{1},t)=\frac{dx_{2}}{dt}
\end{equation}
Lema\^itre propose alors d'appliquer la m\'ethode d'it\'eration rationnelle aux trois fonctions $x_{1}(t)$, $x_{2}(t)$ et $x_{3}(t)$. L'objectif de Georges Lema\^itre est en fait la question de la r\'esolution approch\'ee d'\'equations diff\'erentielles. Il publie dans le m\^eme num\'ero du {\it Bulletin de la Classe Scientifique} un article sur l'application de la m\'ethode d'it\'eration rationnelle \`a l'int\'egration d'une \'equation diff\'erentielle ordinaire \cite{Lem2} dont nous allons parler.

\subsection{Un second article de Georges Lema\^itre}
Goerges Lema\^itre propose donc dans \cite{Lem2}  une m\'ethode d'int\'egration num\'erique pour la r\'esolution d'\'equations diff\'erentielles ordinaires. Il d\'eveloppe son propos en traitant l'exemple de l'\'equation 
\begin{equation}
y'=\frac{dy}{dx}=2y^2(y-x) 
\end{equation}
o\`u $y$ est fonction de $x$. Il s'agit pour lui de calculer des valeurs approch\'ees de la solution $y$ qui tend vers $x$ quand $x$ tend vers $+\infty$. Georges Lema\^itre n'explique pas pourquoi une telle solution existe. Georges Lema\^itre va cependant consid\'erer l'isocline $1$ ($1=2y^2(y-x)$) et la droite $y=x$ (l'isocline $0$). En fait, il y a un anti-entonnoir compris entre l'isocline $1$ et la droite. Le signe de la d\'eriv\'ee partielle de $y^2(y-x)$ par rapport \`a $y$ permet de conclure qu'il y a une solution et une seule pi\'eg\'ee dans cet anti-entonnoir\footnote{Pour toutes ces notions ont peu se r\'ef\'erer \`a \cite{Hubbard}.}. 

Nous avons essay\'e de voir ce qui se passe avec un ordinateur. Maple \footnote{sur la version 14 de Maple}, le logiciel choisi, utilise une m\'ethode de Runge-Kutta-Fehlberg (RKF45)\footnote{qui est une m\'ethode explicite.}. Nous n'avons pas pu mettre en \'evidence cette solution en travaillant dans le sens positif, i.e. sur l'intervalle $\left\lbrack 0,+\infty\right\lbrack$. Cependant, on trouve des valeurs semblables jusqu'\`a $2.2$ \`a celles trouv\'ees par Lema\^itre et r\'esum\'ees dans la tableau se trouvant \`a la fin de l'article (il calcule les valeurs sur l'intervalle $\left\lbrack 0;3\right\rbrack$. La valeur de la condition initiale donn\'ee par Lema\^itre, $x(0)=0.618343$ joue bien un r\^ole de \og s\'eparateur \fg dans le sens o\`u si l'on prend une valeur l\'eg\`erement sup\'erieure, la solution sort de l'anti-entonnoir par le haut, tandis que si l'on prend une valeur inf\'erieure, la solution traverse l'isocline $0$. On peut mettre en \'evidence la solution en travaillant dans le sens n\'egatif, i.e. sur l'intervalle $\left\rbrack -\infty,0\right\rbrack$ car la solution est alors stable dans ce sens\footnote{Nous remercions le professeur Ernst Hairer de nous avoir signal\'e cela et d'avoir fait des essais num\'eriques en utilisant son code RADAU5 (Runge-Kutta implicite). En partant d'une valeur initiale $y(x^{\star}) = x^{\star}$, il obtient dans la direction n\'egative les valeurs suivante pour $y(0)$ :  si $x^{\star}=  1.0$, $y(0 )= 0.602680285037$ ; si $x^{\star}=  3.0$,    $y(0 )= 0.618340077402$ ;  si $x^{\star}= 10^1$,   $y(0 )= 0.618340077404$ ;
si $x^{\star}= 10^2$,   $y(0) = 0.618340077404$ ;
si $x^{\star}= 10^5$ ,  $y(0) = 0.618340077404$ ;
si $x^{\star} = 10^8$,   $y(0) = 0.618340077404$. La  pr\'ecision est  de $3*10^{-15}$ et le nombre de pas d'int\'egration est de $582$ pour $x^{\star}=10^8$.}.

Georges Lema\^itre n'explique pas pourquoi il a choisi cette \'equation\footnote{Georges Lema\^itre s'int\'eresse particuli\`erement \`a cette \'epoque au probl\`eme de St\"ormer et d\'eveloppe des algorithmes pour les machines m\'ecaniques de calcul num\'erique qu'il utilise \`a Louvain.}. On peut imaginer qu'il ait pu choisir cette \'equation parce que les m\'ethodes classiques d'int\'egration num\'erique ne fonctionnent pas de fa\c con satisfaisante. Il applique donc la m\'ethode d\'ecrite \`a la fin de \cite{Lem1} : il part d'une courbe  li\'ee \`a l'\'equation diff\'erentielle qui donne la \og {\it premi\`ere approximation $y_{1}$} \fg, il d\'erive et int\'egre $y_{1}$ pour obtenir deux nouvelles fonctions $y_{0}$ et $y_{2}$. Georges Lema\^itre applique ensuite l'it\'eration rationnelle aux valeurs $y_{0}(x)$, $y_{1}(x)$ et $y_{2}(x)$. Pour choisir la fonction initiale $y_{1}$, Lema\^itre consid\`ere la famille des isoclines $p=2y^2(y-x)$ dont il forme l'\'equation diff\'erentielle associ\'ee :
\begin{equation}
y'=\frac{dy}{dx}=\frac{y}{3y-2x}.
\end{equation}
Par des conditions g\'eom\'etriques (qualitatives), Lema\^itre consid\`ere la \og {\it courbe limite} \fg donn\'ee par la condition :
\begin{equation}
y'=\frac{dy}{dx}=\frac{y}{3y-2x}=p=2y^2(y-x)
\end{equation}
soit
\begin{equation}
2y_{1}(y_{1}-x)(3y_{1}-2x)=1
\end{equation}
courbe somme toute proche de l'isocline $1$. Il lui reste \`a la param\'etrer.

La m\'ethode de Georges Lema\^itre est dont un m\'elange de m\'ethodes qualitatives (choix de $y_{1}$), de la m\'ethode it\'erative de Picard (avec aussi des m\'ethodes pour les calculs annexes) et du $\Delta^{2}$ d'Aitken pour am\'eliorer les r\'esultats.

\section{Conclusion}
L'objectif de cette note n'\'etait pas de r\'egler une question de primaut\'e d'un r\'esultat scientifique. Les exemples abondent en math\'ematiques de m\'ethodes red\'ecouvertes puis am\'elior\'ees ou de r\'esultats prouv\'es ind\'ependamment. Les pr\'eoccupations des auteurs et les dynamiques de recherches ont leur importance dans ce ph\'enom\`ene. Dans le cas de l'it\'eration 
rationnelle, c'est la r\'esolution de certaines \'equations diff\'erentielles qui incite Georges Lema\^itre \`a chercher un moyen d'acc\'elerer une m\'ethode connue (celle de Picard) et qui soit facile \`a mettre en oeuvre, en un mot \`a impl\'ementer sur une machine \`a calculs m\'ecanique.

Les motivations des trois math\'ematiciens \'etaient bien distinctes. Aitken se d\'emarque, son objectif de d\'epart est le calcul approch\'e des racines d'une \'equation polyn\^omiale \`a une variable. Steffensen et Lema\^itre se rejoignent : il s'agit d'acc\'el\'erer la convergence d'un processus it\'eratif. Mais Steffensen en reste au stade des it\'er\'ees d'une fonction \`a une variable, tandis que Lema\^itre applique sa m\'ethode \`a des fonctions de plusieurs variables, mais surtout \`a des \'equations diff\'erentielles

L'isolement scientifique et les affres de la guerre expliquent en partie que Lema\^itre n'ait pas eu connaissance de l'existence des travaux de Aitken et de Steffensen. Ajoutons que Lema\^itre avait une fa\c con bien \`a lui de travailler, en citant Lambert\footnote{Lambert op.cit., page 147} :
\begin{quotation}
\begin{it}
Remarquons aussi un d\'etail qui est caract\'eristique de la mani\`ere dont Lema\^itre travaillait. Il poursuit les id\'ees qui l'int\'eressent en ne se souciant gu\`ere de la litt\'erature existante ou des groupes qui pourraient aborder les m\^emes sujets que lui.
\end{it}
\end{quotation}

\newpage
\section{Annexes}
\subsection{Compte-rendu du Jahrbuch de l'article de Aitken}
\begin{quotation}
\begin{it}
%Jahrbuch
Aitken, A. C.

On Bernoulli's numerical solution of algebraic equations, Proceedings Royal Soc. Edinburgh 46, 289-305.

Die {\it Bernoulli}sche Methode zur numerischen Aufl\"osung einer
    algebraischen Gleichung 
$$ g (x) = a_0 x^n + a_1 x^{n-1} + \cdots
    + a_n = 0 
$$ 
beruht auf folgendem Satze: Wird bei beliebigen
    Anfangsgliedern durch die Rekursionsgleichung 
$$ a_0 f(t + n)
    + a_1f(t + n - 1) + \cdots + a_n f(t) = 0 
$$ 
die Folge $f(t)$
    bestimmt, so ist 
$$ \lim_{t\to\infty} \frac{\overline{f(t+1)}}{f(t)}
    = z_1, 
$$
 wobei $z_1$ die Wurzel vom gr\"o{\ss}ten absoluten
    Betrage bedeutet; allerdings mu{\ss} noch vorausgesetzt werden,
    da{\ss} die andern Wurzeln s\"amtlich von kleinerem Betrage sind
    als $z_1$. \par Diese {\it Bernoulli}sche Formel verallgemeinert
    nun Verf. dahingehend, da{\ss} die M\"oglichkeit gegeben wird,
    s\"amtliche Wurzeln gleichzeitig zu bestimmen. \par Setzt man
    n\"amlich $f_m(t)$ gleich, der Determinante 
$$
 | f(t + \lambda
    - \kappa)| \qquad (\kappa, \lambda = 1, 2, \dots, m), $$
    so ist $$ \lim_{t\to\infty}\frac{f_m(t+1)}{f_m(t)} = z_1 z_2
    \dots z_m 
$$
 wobei $z_1, z_2, \dots, z_n$ die Wurzeln von $g(x)
    = 0$ in absteigender Ordnung nach ihrem absoluten Betrage bedeuten
    und vorausgesetzt wird, da{\ss}$|z_m| > |z_{m+1}|$ sei. \par
    Die Praxis bietet auch f\"ur die Falle, in denen die allgemeine
    Formel versagt, M\"oglichkeiten, die Wurzeln einer algebraischen
    Gleichung nach der Regel zu bestimmen, die durch den angegebenen
    Satz angedeutet ist. Diese M\"oglichkeiten werden in den weiteren
    Abschnitten der Arbeit an Hand einiger Beispiele diskutiert.
    \par Im weiteren Verlauf der Untersuchungen werden aus dem Vorangehenden
    noch Methoden hergeleitet, sch\"arfere Approximationen der Wurzeln
    zu erzielen. \par Zum Schlu{\ss} wird die Methode des Verf. noch
    mit \"andern verglichen, insbesondere mit der Methode von {\it
    F\"urstenau} (Darstellung der reellen Wurzeln algebraischer Gleichungen
    durch Determinanten der Coefficienten, Marburg 1860), die als
    ein Spezialfall (bei passender Wahl der Anfangselemente der Folge
    $f (t)$) erscheint. (IV 17.)

W.; Dr. (K\"onigsberg i. Pr.)
\end{it}
\end{quotation}
\subsection{Compte-rendu du Jahrbuch de l'article de Steffensen}
%Jahrbuch
\begin{quotation}
\begin{it}
Steffensen, J. F.\\
Remark on iteration, Skand. Aktuarietidskrift 16, 64-72,1933

Wird die zu l\"osende Gleichung auf die Form $x = f(x)$ gebraucht,
    und setzt man $x_{\nu + 1} = f(x_{\nu })$, so ergibt sich ein
    guter Ann\"aherungswert aus einer Folge dreier Werte $x_0, x_1,
    x_2$ durch $$x = x_0 - \frac {(\Delta x_0)^2}{\Delta ^2x_0}.$$

Lorey, W.; Prof. (Frankfurt am Main)
\end{it}
\end{quotation}
\subsection{Compte-rendu du Jahrbuch de l'article de Lema\^itre sur l'it\'eration rationnelle}
%Jahrbuch
\begin{quotation}
\begin{it}
Anschlie{\ss}end an eine von Gau{\ss} in der Theoria motus entwickelte
    Methode gibt Verf. ein Verfahren an, aus mehreren durch einfache
    Iteration gewonnenen N\"aherungswerten einen neuen besseren Wert
    zu gewinnen. F\"ur eine Gleichung mit einer Unbekannten ergibt
    sich unter der Annahme, da{\ss} die Funktion in dem in Frage
    kommenden Gebiet linear angen\"ahert werden darf, aus der Vertr\"aglichkeitsbedingung
    der drei Gleichungen $$ x_2 = f(x_1) = a + bx_1, \quad x_3 =
    f(x_2) = a + bx_2, \quad x = f(x) = a + bx $$ der neue N\"aherungswert
    $$ x = \frac{x_1 x_3 - x_2^2}{x_1 - 2x_2 + x_3}. $$ Der Vorteil
    dieser Formel ist der, da{\ss} man auch zum Ziel kommt, wenn
    die Folge der N\"aherungswerte $x_1$, $x_2$, $x_3$ nicht gegen
    $x$ konvergiert (das Verfahren stimmt mit einem von {\it Steffensen},
    Skand. Aktuarietidskr. 16 (1933), 64-72 (F. d. M. 58,
    535) gegebenen \"uberein). Bei Gleichungen mit mehreren Unbekannten
    erh\"alt man f\"ur jede der Unbekannten eine entsprechende Gleichung.
    \par Zum Schlu{\ss} wird kurz auf die Anwendung des Verfahrens
    auf Differentialgleichungen eingegangen. Dabei unterscheidet
    Verf. eine Iteration durch Integration und eine durch Differentiation.

Willers, F.; Prof. (Dresden)
\end{it}
\end{quotation}

\subsection{Compte-rendu du Jahrbuch de l'article de Lema\^itre sur la r\'esolution d'une \'equation diff\'erentielle par l'it\'eration rationnelle}
Int\'egration d'une \'equation diff\'erentielle par it\'eration rationnelle.

Acad. Belgique, Bull. Cl. Sci. (5) 28, 815-825.

Das in der vorstehend besprochenen Arbeit beschriebene Verfahren
    wird auf die Gleichung $$ \frac{dy}{dx} = 2y^2 (y - x) $$ angewendet.
    Aus einer ersten N\"aherung $y_1$ werden zwei weitere durch $$
    \frac{1}{y_2} = 2 \int_{x}^{\infty} (y_1 - x)dx \ \mbox{
    und } \ y_0 - x = -\frac{1}{2}\frac{d}{dx}\frac{1}{y_1}
    $$ gewonnen und aus diesen in der angegebenen Art eine neue N\"aherung
    berechnet. Gesucht wird eine L\"osung im Intervall $0 \leq x
    \leq \infty$, die die Gerade $y = x$ als Asymptote hat. Ausgegangen
    wird von der Wendepunktskurve der Differentialkurven $$ 2y_1(y_1
    - x)(3y_1 - 2x) = 1, $$ die zun\"achst punktweise mit dem Argumentsschritt
    $\Delta x =$ 0,1 von $x = 0$ bis $x = 3$ berechnet wird. Dazu
    wird $y_0$ gebildet und berechnet und weiter durch numerische
    Integration $y_2$. Die aus diesen drei N\"aherungen berechnete
    neue ist auf sechs Dezimalen genau. F\"ur gro{\ss}e $x$ wird
    nach der Methode der unbestimmten Koeffizienten eine asymptotische
    Reihe bestimmt, aus der sich f\"ur $x = 3$ bis auf sechs Stellen
    der gleiche Wert ergibt wie der oben berechnete. Angabe der Resultate
    in Tabellen.

Willers, F.; Prof. (Dresden)

\end{document}